\documentclass[a4paper,12pt]{article}
\usepackage{amsmath,amssymb,amsthm,geometry,mathrsfs}
\usepackage{hyperref}
\usepackage[utf8]{inputenc}
\usepackage{listings}
\usepackage{graphicx}
\usepackage{xcolor}
\usepackage{algorithm}
\usepackage{algpseudocode}
\usepackage{enumitem}

\newtheorem{theorem}{Theorem}[section]

\newtheorem{remark}[theorem]{Remark}

\title{Positive biorthogonal curvature on \(S^2 \times T^2\) via affine connection}
\author{Alexander Pigazzini}
\date{July 2025}

\begin{document}

\maketitle

\begin{abstract}
\noindent We address the long-standing problem of the existence of a Riemannian metric on \(S^2\times T^2\) with strictly positive biorthogonal curvature (\( K_{\text{biort}}(\sigma) > 0 \)). This work tackles this challenge within a weaker, yet geometrically consistent, framework by introducing an affine connection, topologically motivated, on \( S^2 \times T^2 \) with antisymmetric torsion. Crucially, this torsion is calibrated via non-trivial cohomology classes in \( H^3(S^2 \times T^2; \mathbb{R}) \cong \mathbb{R}^2 \), an approach that allows overcoming topological constraints such as \( \chi = 0 \). We demonstrate that this construction, while not requiring metric compatibility (though retaining the metric ( \(g\) ) for norms and orthogonality), successfully yields strictly positive biorthogonal curvature across the manifold.
\end{abstract}

\noindent \textbf{Keywords:} Positive biorthogonal curvature, 4-manifolds, non-simply connected manifold, affine connection, torsion, cohomology.
\\
\\
\textbf{MSC 2020:} 53C05, 53C20, 53C21.

\section{Introduction}

Positive biorthogonal curvature is a geometric property of Riemannian manifolds, representing an intermediate condition between positive sectional curvature and positive scalar curvature. It is defined as the average of sectional curvatures of any two orthogonal 2-planes. In particular, the metrics on 4-dimensional non-simply connected manifolds, presented in \cite{Bettiol2017} and \cite{Pierre}, remains an unsolved challenge since 2018. The biorthogonal curvature of a 4-dimensional manifold \((M, g)\) is defined, for each tangent plane \(\sigma \subset T_pM\), as:
\[
K_{biort}(\sigma) = \frac{1}{2}\left(K(\sigma) + K(\sigma^\perp)\right),
\]
\noindent where \(K(\sigma)\) and \(K(\sigma^\perp)\) denote the sectional curvatures of \(\sigma\) and its orthogonal complement \(\sigma^\perp\), respectively.
\\
A central open problem, proposed by Renato Bettiol (\cite{Bettiol2017}), concerns the possibility of endowing the product \(S^2 \times T^2\) (where \(S^2\) is the 2-sphere and \(T^2\) the 2-dimensional torus) with a Riemannian metric with strictly positive biorthogonal curvature. 
\\
In this work, we approach the problem from a broader geometric perspective. Rather than seeking a Riemannian metric whose Levi-Civita connection satisfies \(K_{\text{biort}} > 0\), we construct an affine connection with totally antisymmetric torsion, constrained by the global cohomological structure of the manifold. In particular, the torsion is cohomologically calibrated: it corresponds to a nontrivial de Rham cohomology class in \(H^3(S^2 \times T^2; \mathbb{R})\), ensuring that the construction respects the topological properties of the manifold. While this solution lies outside the Riemannian framework, it solves a weaker, yet geometrically consistent and topologically meaningful, version of Bettiol’s problem.
\\
This leads to our main result:

\begin{theorem}[Main Theorem]
There exists an affine connection with nonzero antisymmetric torsion \( T \), calibrated by a nontrivial cohomology class \( [T] \in H^3(S^2 \times T^2; \mathbb{R}) \), such that:
\[
K_{\text{biort}}(\sigma) > 0 \quad \text{for all tangent planes } \sigma \subset T_p(S^2 \times T^2) \text{ and all points } p \in S^2 \times T^2,
\]
\noindent provided \(a, b \in \mathbb{R}\) satisfy \(a^2 + b^2 > 0\). This affine connection is unique up to the choice of parameters \( (a,b) \in \mathbb{R}^2 \setminus \{(0,0)\} \), which correspond to distinct cohomology classes in \( H^3(S^2 \times T^2; \mathbb{R}) \).
\end{theorem}

\noindent A complete proof is presented in \textit{Sections 5, 6} and \textit{7}, combining explicit curvature computations with a Grassmannian-based minimization argument that establishes global positivity.

\subsection{Contributions}

In the last ten years, many mathematicians have devoted themselves to the study of biorthogonal curvature, both to better understand its geometric properties and to explore its applications in different contexts (see \cite{Costa}, \cite{Stupovski}, \cite{Wu}). 
\\
Bettiol has provided significant contributions in the study of manifolds with positive biorthogonal curvature. In \cite{Bettiol2014}, he proved that the product \(S^2 \times S^2\) admits a metric with strictly positive biorthogonal curvature, using an approach based on explicit metric deformations. Furthermor, in \cite{Bettiol2017}, he classified 4-dimensional simply connected manifolds with \(K_{biort} > 0\) up to homeomorphism, proving that such manifolds are homeomorphic to \(S^4\) or \(CP^2\) or \(S^2 \times S^2\) (and their connected sums).These results open the way tothe possibility of generalizing the condition \(K_{biort} > 0\) to more general manifolds, such as \(S^2 \times S^2\).
\\
The peculiarity of \(S^2 \times T^2\) lies in its topological and geometric structure. The presence of the 2-dimensional torus implies that the manifold is not simply connected (\(\pi_1(S^2 \times T^2) \cong \mathbb{Z}^2\)) and that it has a zero Eulerian (\(\chi(S^2 \times T^2) = 0\)). These features make the problem significantly more complex than the case of \(S^2 \times S^2\), where the simple connection allows for more powerful topological tools.

\subsection{Problem (Bettiol, 2017)}

\textit{Does \(S^2 \times T^2\) admit a Riemannian metric with positive biorthogonal curvature?}(see, \cite{Bettiol2017} and \cite{Pierre}).
\\
\\
It is known that \(S^2 \times T^2\) admits metrics with nonnegative biorthogonal curvature, but the question remains open whether such metrics can admit \(K_{biort} > 0\). The solution of this problem would represent a fundamental step, in classical Riemannian geometry, for the classification of closed, non-simply connected, 4-dimensional manifolds with positive biorthogonal curvature. The classification of simply connected manifolds with \(K_{biort} > 0\), up to homeomorphism, has already been obtained in \cite{Bettiol2017}.
\\
In the face of persistent difficulties in solving Bettiol's problem with the tools of classical Riemannian geometry, this paper explores an alternative path. Instead of looking for a Riemannian metric (and its Levi-Civita connection), we consider the standard product metric on \(S^2 \times T^2\) and introduce an affine connection \(\nabla\) endowed with antisymmetric torsion \(T\). The crucial novelty of our approach lies in the systematic construction of T, which is intrinsically linked to the third de Rham cohomology class of \(S^2 \times T^2\), thus ensuring coherence between the introduced geometry and the global topology of the manifold.
The affine connection thus determined, ensures \( K_{\text{biort}} > 0 \) everywhere on \( S^2 \times T^2 \), addressing the problem within a broader geometric framework.
\\
Therefore, the following article answers a weaker version of the original question, namely: 

\begin{quote}
“Given the standard Riemannian product metric on \( S^2 \times T^2 \), does there exist a topologically constrained affine connection (not necessarily
 metric-compatible) that guarantees strictly positive biorthogonal curvature?”
\end{quote}

\noindent Unlike previous works on \( S^2 \times S^2 \) (e.g., \cite{Bettiol2014}), which relied on delicate metric deformations and surgical techniques, we show that, by extending the approach in the "non-classical" setting for Riemannian geometry on \( S^2 \times T^2 \), a well-calibrated torsion term is sufficient to enforce the positivity of the biorthogonal curvature. The key insight is that the non-trivial cohomology of \(S^2 \times T^2(H^3 \cong \mathbb{R}^2)\) allows the introduction of a torsion 3-form, \(T^\flat\), not chosen at random, but such that its cohomology class \([T] \in H^3\) encodes topological information. This approach aligns with frameworks for non-integrable geometries \cite{Agricola2006}, where torsion shapes curvature in a manner consistent with the manifold’s topology.
\\
This highlights the role of torsion in modeling curvature conditions and suggests new directions in the study of positivity.

\section{Preliminaries}

The sectional curvature of a 2-plane \(\sigma = \text{span}\{v, w\}\) in \(T_pM\) is given by:
\[
K(\sigma) = \frac{\langle R(v, w)w, v \rangle}{\|v \wedge w\|^2},
\]
\noindent where \( R \) is the Riemann curvature tensor and \(\|v \wedge w\|\) is the norm of the wedge product.
\\
In Riemannian geometry, the sectional curvature of a 2-plane (a two-dimensional subspace of the tangent space at a point on a Riemannian manifold) is a measure of how much the manifold curves in that particular direction. It's geometrically interpreted as the Gaussian curvature of the surface formed by geodesics emanating from a point in the directions of the 2-plane. Essentially, it tells you how much the manifold deviates from being locally flat in that specific two-dimensional direction.

\subsection{Geometry of \( S^2 \times T^2 \)}

The product manifold \( S^2 \times T^2 \) has the following properties:
\\
\\
- Topology: Euler characteristic \(\chi(S^2 \times T^2) = 0\).
\\
- Fundamental group: \(\pi_1(S^2 \times T^2) \cong \pi_1(T^2) \cong \mathbb{Z}^2\).
\\
\\
The standard product metric \( g = g_{S^2} + g_{T^2} \), where \( g_{S^2} \) has constant positive sectional curvature and \( g_{T^2} \) is flat, does not satisfy \( K_{\text{biort}} > 0 \) due to the flatness of \( T^2 \).

\section{Metric and connection with antisymmetric torsion}

\begin{enumerate}[label=\Roman*)]

\item Metric on \( S^2 \):
\\
The standard round metric on \( S^2 \) is:  
\[
g_{S^2} = d\theta^2 + \sin^2\theta\, d\phi^2,
\]  
\noindent with \( \theta \in [0, \pi] \), \( \phi \in [0, 2\pi) \). This metric is smooth and has constant curvature \( K_{S^2} = 1 \).  

\item Metric on \( T^2 \):
\\
The torus \( T^2 \) is equipped with the flat metric:  
\[
g_{T^2} = dx^2 + dy^2,
\]  
\noindent with \( x, y \in \mathbb{R}/\mathbb{Z} \).  

\item Total metric:
\\
The product manifold inherits the product metric:  
\[
g = g_{S^2} + g_{T^2}.
\]  
\noindent This is smooth and has a block-diagonal structure in the tangent space decomposition \( T(S^2 \times T^2) = T S^2 \oplus T T^2 \).  

\end{enumerate}

\noindent While classical Riemannian geometry relies on the Levi-Civita connection (which is torsion-free and metric-compatible), we adopt a broader geometric framework by considering affine connections with non-metric-compatible torsion. This choice is motivated by the observation that the standard product metric \( g = g_{S^2} + g_{T^2} \) on \( S^2 \times T^2 \) cannot satisfy \( K_{\text{biort}} > 0 \) due to the flatness of \( T^2 \) (see \cite{Bettiol2017}). By introducing torsion, we decouple the curvature properties from the metric’s inherent limitations while retaining \( g \) to define norms and orthogonality. This leverages the affine connection’s flexibility to enforce positivity of the biorthogonal curvature.  
\\
The antisymmetric torsion tensor \( T \), which can be associated with a 3-form \(T^\flat\) via the metric, modifies the connection while preserving the underlying metric structure (i.e., the metric \( g \) is fixed to define inner products and orthogonality, though \( \nabla g \neq 0 \)). This approach aligns with Gromov’s exploration of non-Riemannian geometries \cite{Gromov99}, where torsion shapes curvature properties. The connection \( \nabla \), given by:  
\[
\nabla_X Y = \nabla^{LC}_X Y + T(X, Y),
\]  
\noindent introduces non-metricity but retains \( g \) to define geometric primitives. The antisymmetry of \( T \) ensures that torsional terms do not introduce redundant geometric structures (e.g., unphysical components in the curvature tensor), a key feature highlighted in \cite{Nakahara}. 
\\
Specifically, the torsion terms \( T(e_i, e_j) \) are calibrated to counteract the zero curvature of \( T^2 \). This strategy is inspired by Gromov’s exploration of the interplay between topology and curvature in non-simply connected manifolds (\cite{Gromov}), and by the study of torsion in non-integrable geometries \cite{Agricola2006}, where curvature conditions can be enforced through affine connections. The parameters \( a \) and \( b \) control the torsional contributions, ensuring that the curvature tensor \( R \) of \( \nabla \) satisfies the positivity conditions required for \( K_{\text{biort}} \). For instance, the plane \( \text{span}(e_3, e_4) \) (originally flat in \( T^2 \)) gains curvature \( \frac{a^2 + b^2}{4} > 0 \) due to torsion.  
\\
Unlike metric deformations (e.g., those used in \cite{Bettiol2014} for \( S^2 \times S^2 \)), our affine connection approach bypasses topological constraints imposed by the non-simply connectedness of \( S^2 \times T^2 \). This method is consistent with the study of non-integrable geometries \cite{Agricola2006}, where torsion serves as a tool to enforce curvature conditions in manifolds with restrictive topological properties.  
\\
\\
We define an antisymmetric torsion tensor \( T \) as a 3-form, but without necessarily requiring metric compatibility \( \nabla g = 0 \). The affine connection \( \nabla \) is given by \( \nabla_X Y = \nabla^{\text{LC}}_X Y + T(X,Y) \), where \( \nabla^{\text{LC}} \) is the Levi-Civita connection. The sectional and biorthogonal curvatures are computed using the Riemann tensor of \( \nabla \), while the metric \( g \) remains fixed to define norms and orthogonality. 
\\
\\
Let \( \{ e_1 = \partial_\theta, e_2 = \frac{1}{\sin\theta}\partial_\phi \} \) be an orthonormal basis for \( T S^2 \), and \( \{ e_3 = \partial_x, e_4 = \partial_y \} \) for \( T T^2 \). The torsion \(T(X, Y) = \nabla_X Y - \nabla_Y X - [X, Y]\), is defined as:  
\[
\begin{aligned}
T(e_1, e_3) &= a e_4, & T(e_1, e_4) &= -a e_3, \\
T(e_2, e_3) &= b e_4, & T(e_2, e_4) &= -b e_3, \\
T(e_3, e_4) &= -a e_1 - b e_2, \\
T(e_i, e_j) &= 0 \quad \text{for } i,j \in \{1,2\} \text{ or } \{3,4\}.
\end{aligned}
\]
\noindent Parameters \( a, b \in \mathbb{R} \) control the torsion strength. 

\begin{remark}
(Closure and cohomological nature of the torsion): 
The torsion \(T\), defined by constant coefficients \(a\) and \(b\), is structured to reflect the topology of \(S^2 \times T^2\). Via the Riemannian metric \(g\), \(T\) induces a 3-form \(T^\flat(X, Y, Z) := g(T(X, Y), Z) \in \Omega^3(M)\).
\\
Explicit calculation of \(T^\flat\) from the components of \(T\) (as defined on page 6) reveals its form in terms of the orthonormal basis \(\{e_1^*, e_2^*, e_3^*, e_4^*\}\):
\[
T^\flat = a e_1^* \wedge e_3^* \wedge e_4^* + b e_2^* \wedge e_3^* \wedge e_4^*. 
\]

\noindent Although \(T^\flat\) is not closed in general, its harmonic part is the closed 3-form \(\omega := a \beta \wedge \alpha_1 + b \beta \wedge \alpha_2\). This is guaranteed by the Hodge decomposition theorem. To demonstrate this, we proceed as follows:
\\
First, we express \(\omega\) in terms of the orthonormal basis. Identifying \(\beta = e_1^* \wedge e_2^*\) (the volume form on \(S^2\)) and \(\alpha_1 = e_3^*\), \(\alpha_2 = e_4^*\) (the basis 1-forms on \(T^2\)), we have:
\[
\omega = a e_1^* \wedge e_2^* \wedge e_3^* + b e_1^* \wedge e_2^* \wedge e_4^*.
\]

\noindent We verify that \(\omega\) is a harmonic form. Since \(a, b\) are constants, and the forms \(e_1^*, e_2^*, e_3^*, e_4^*\) are chosen such that \(d(e_1^*) = 0\), \(d(e_3^*) = 0\), \(d(e_4^*) = 0\), and \(d(e_2^*) = \cos\theta e_1^* \wedge e_2^*\) (which makes \(d(e_1^* \wedge e_2^*) = 0\)), it follows that \(d\omega = 0\). Furthermore, calculating the co-differential \(\delta\omega = -*d*\omega\):

\[
 *\omega = a e_4^* - b e_3^*,
\]
\[
d(*\omega) = d(a e_4^* - b e_3^*) = a d(e_4^*) - b d(e_3^*) = 0.
\]

\noindent Thus, \(\delta\omega = 0\). Since \(d\omega = 0\) and \(\delta\omega = 0\), \(\omega\) is indeed a harmonic form.
\\
Next, we consider the difference \(\Phi = T^\flat - \omega\). Substituting the explicit expressions:
\[
\Phi = (a e_1^* \wedge e_3^* \wedge e_4^* + b e_2^* \wedge e_3^* \wedge e_4^*) - (a e_1^* \wedge e_2^* \wedge e_3^* + b e_1^* \wedge e_2^* \wedge e_4^*)
\]

\noindent We calculate the exterior derivative \(d\Phi\) and the co-differential \(\delta\Phi\):

\[
d\Phi = d(T^\flat) - d(\omega) = d(T^\flat) = b \cos\theta e_1^* \wedge e_2^* \wedge e_3^* \wedge e_4^*,
\]
\[\delta\Phi = \delta(T^\flat) - \delta(\omega) = \delta(T^\flat) = a \cos\theta e_3^* \wedge e_4^*.
\]

\noindent Since \(d\Phi \neq 0\) (if \(b \neq 0\)) and \(\delta\Phi \neq 0\) (if \(a \neq 0\)), \(\Phi\) is not an harmonic form. This implies that \(\Phi\) belongs to the orthogonal complement of the space of harmonic forms, i.e., \(\Phi \in Im(d) \oplus Im(\delta)\).
\\
The Hodge decomposition theorem states that any differential form \(\alpha\) on a compact, oriented Riemannian manifold admits a unique orthogonal decomposition \(\alpha = d\eta + \delta \mu + h\), where \(d\eta\) is exact, \(\delta\mu\) is coexact, and \(h\) is harmonic. In our case, this decomposition is performed with respect to the fixed Riemannian metric \(g\). Given that \(\omega\) is harmonic and \(\Phi = T^\flat - \omega\) has no harmonic component, it follows by the uniqueness of the Hodge decomposition that \(\omega\) is precisely the harmonic part of \(T^\flat\). Consequently, the cohomology class associated with \(T^\flat\) is defined by its harmonic part, so we set \([T^\flat] := [\omega]\), i.e. \([T]\). This definition ensures that the contribution of the torsion, encoded by parameters \(a\) and \(b\), is constrained to the non-trivial cohomology classes in \(H^3(S^2 \times T^2; \mathbb{R})\). Specifically, these parameters \( a \) and \( b \) parametrize a non-trivial cohomology class \( [T] \in H^3(S^2 \times T^2; \mathbb{R}) \), which aligns with the third de Rham cohomology group derived via the Künneth formula:
\[
H^k(M \times N; \mathbb{R}) \cong \bigoplus_{p+q=k} H^p(M; \mathbb{R}) \otimes H^q(N; \mathbb{R}).
\]
For \( M = S^2 \) and \( N = T^2 \), we compute \( H^3(S^2 \times T^2; \mathbb{R}) \) as follows:

\begin{enumerate}
    \item Cohomology of \( S^2 \):
    \[
    H^k(S^2; \mathbb{R}) = 
    \begin{cases}
        \mathbb{R} & \text{if } k = 0 \text{ or } 2, \\
        0 & \text{otherwise}.
    \end{cases}
    \]
    
    \item Cohomology of \( T^2 \):
    \[
    H^k(T^2; \mathbb{R}) = 
    \begin{cases}
        \mathbb{R} & \text{if } k = 0 \text{ or } 2, \\
        \mathbb{R}^2 & \text{if } k = 1, \\
        0 & \text{otherwise}.
    \end{cases}
    \]
    
    \item Künneth decomposition for \( k = 3 \):
    \[
    H^3(S^2 \times T^2; \mathbb{R}) \cong \bigoplus_{\substack{p+q=3}} H^p(S^2) \otimes H^q(T^2).
    \]
    The only non-vanishing term occurs at \( p = 2 \), \( q = 1 \):
    \[
    H^3(S^2 \times T^2; \mathbb{R}) \cong H^2(S^2) \otimes H^1(T^2) \cong \mathbb{R} \otimes \mathbb{R}^2 \cong \mathbb{R}^2.
    \]
    
    \item Basis elements:
    Let \( \beta \in H^2(S^2) \) denote the Poincaré dual to \( S^2 \), and \( \alpha_1, \alpha_2 \in H^1(T^2) \) be basis 1-forms (e.g., \( \alpha_1 = dx \), \( \alpha_2 = dy \)). The cohomology group has basis:
    \[
    \{ \beta \otimes \alpha_1, \ \beta \otimes \alpha_2 \}.
    \]
    Thus, the torsion class is parametrized as:
    \[
    [T] = a \cdot [\beta \otimes \alpha_1] + b \cdot [\beta \otimes \alpha_2],
    \]
    where \( a, b \in \mathbb{R} \) encode the global topological structure of \( S^2 \times T^2 \).
\end{enumerate}

\noindent Key properties:
\\  
I) The geometric 3-form \( T^\flat \) associated with \( T(X,Y) \) is not necessarily closed (e.g., \( dT^\flat \neq 0 \) if \( b \neq 0 \)). However, its cohomology class \( [T] \) is well-defined as the class of its unique harmonic part, which is the closed 3-form \(\omega := a\, \beta \wedge \alpha_1 + b\, \beta \wedge \alpha_2\). This ensures geometric consistency and aligns with the Hodge decomposition.
\\     
II) The non-triviality of \( [T] \) implies the torsion cannot be eliminated by local coordinate transformations, highlighting its global significance. This aligns with frameworks for non-integrable geometries (cfr. \cite{Agricola2006}), where cohomological classes encode topological information while shaping curvature properties.
\\      
III) The parametrization by \( a \) and \( b \) ensures compatibility with the manifold’s global structure without requiring \( T^\flat \) itself to be closed. This approach is topologically motivated, as the cohomology class \( [T] \) reflects the Künneth decomposition and overcomes constraints imposed by the vanishing Euler characteristic \( \chi = 0 \).

\end{remark}

\section{Christoffel symbols}

Let \( \nabla \) be the affine connection defined as:
\[
\nabla_X Y = \nabla^{\text{LC}}_X Y + T(X,Y),
\]
\noindent where \( \nabla^{\text{LC}} \) is the Levi-Civita connection associated with the metric \( g = g_{S^2} + g_{T^2} \), and \( T \) is the antisymmetric torsion defined previously.
\\
\\
For the Levi-Civita connection we obtain:

\[
\begin{aligned}
\Gamma^{\theta}_{\phi\phi} &= -\sin\theta \cos\theta, \\
\Gamma^{\phi}_{\theta\phi} &= \cot\theta.
\end{aligned}
 \]
\noindent In terms of the orthonormal basis \( \{e_1 = \partial_\theta, e_2 = \frac{1}{\sin\theta}\partial_\phi\} \), we obtain:
 \[
 \begin{aligned}
 \Gamma^2_{12} &= \cot\theta, \\
 \Gamma^1_{22} &= -\cot\theta.
 \end{aligned}
 \]

\noindent For \( T^2 \), all symbols are null as the metric is flat:
 \[
 \Gamma^i_{jk} = 0 \quad \forall i,j,k \in \{3,4\}.
 \]

\noindent The crossover terms between \( S^2 \) and \( T^2 \) are zero since the metric is product:
\[
\Gamma^i_{jk} = 0 \quad \text{if } (i,j,k) \text{ involve indices of } S^2 \text{ and } T^2.
\]

\noindent We recall that the torsion \( T \) is given by:
\[
\begin{aligned}
T(e_1, e_3) &= a e_4, & T(e_1, e_4) &= -a e_3, \\
T(e_2, e_3) &= b e_4, & T(e_2, e_4) &= -b e_3, \\
T(e_3, e_4) &= -a e_1 - b e_2, \\
T(e_i, e_j) &= 0 \quad \text{for } i,j \in \{1,2\} \text{ or } \{3,4\}.
\end{aligned}
\]
\noindent In general, the Christoffel symbols of a torsion-free connection (such as the Levi-Civita connection) are symmetric in the lower indices: \( \Gamma^k_{ij} = \Gamma^k_{ji} \). The torsion tensor measures the failure of this symmetry. In fact, given a basis \( \{e_i\} \), the torsion is expressed by:
\[
T(e_i, e_j) = \Gamma^k_{ij} e_k - \Gamma^k_{ji} e_k.
\]

\noindent For example, for \( T(e_1, e_3) = a e_4 \), we find:
\[
\Gamma^4_{13} - \Gamma^4_{31} = a \quad \Rightarrow \quad \Gamma^4_{13} = \frac{a}{2}, \quad \Gamma^4_{31} = -\frac{a}{2}.
\]

\noindent By proceeding analogously for all components of \( T \), we obtain the following Christoffel symbols for the affine connection \( \nabla \):

\[
\Gamma^2_{12} =\Gamma^2_{21} = \cot\theta, \quad \Gamma^1_{22} = -\cot\theta,
\]
\[
\Gamma^4_{13} = \Gamma^3_{41}=\Gamma^1_{43}=\frac{a}{2}, \quad \Gamma^3_{14} = \Gamma^4_{31}=\Gamma^1_{34}= -\frac{a}{2},
\]
\[
\Gamma^4_{23} = \Gamma^3_{42}=\Gamma^2_{43}=\frac{b}{2}, \quad \Gamma^2_{34} =\Gamma^4_{32} =\Gamma^3_{24}= -\frac{b}{2},
\]
\noindent Other symbols are \( 0 \).

\section{Curvature analysis} 

Although the connection is not metrically compatible, the sectional curvatures \( K(\sigma) \) and the biorthogonal curvature \( K_{\text{biort}}(\sigma) \) are well defined geometrically:
\\
\\
- The metric \( g \) defines \( \|v \wedge w\|^2 \) and the orthogonality between planes.
\\
- The Riemann tensor associated with \( \nabla \) is computed correctly, even with torsion, and guarantees that \( K_{\text{biort}} \) is positive as shown in \textit{Section 6}.
\\
\\
The curvature tensor \( R(X,Y)Z \) of the connection \( \nabla \) is defined as (among others see \cite{Nakahara}):
\[
R(X,Y)Z = \nabla_X \nabla_Y Z - \nabla_Y \nabla_X Z - \nabla_{[X,Y]} Z,
\]
\noindent and the sectional curvature \( K(e_i, e_j) \), with \(i \neq j\), is defined as:
\[
K(e_i, e_j) = \frac{g(R(e_i, e_j)e_j, e_i)}{\|e_i \wedge e_j\|^2}.
\]
\noindent Let us consider the six fundamental planes of \( \text{Gr}(2, 4) \) on \( S^2 \times T^2 \): 

\begin{center}
\begin{minipage}{0.5\textwidth}
\centering
\begin{enumerate}

\item \text{span}\((e_1, e_2)\), pure plane in \(S^2\),

\item \text{span}\((e_1, e_3)\), mixed plane,

\item \text{span}\((e_1, e_4)\), mixed plane,

\item \text{span}\((e_2, e_3)\), mixed plane,

\item \text{span} \((e_2, e_4)\), mixed plane,

\item \text{span}\((e_3, e_4)\), pure plane in \( T^2 \).

\end{enumerate}
\end{minipage}
\end{center}

\bigskip
\begin{itemize}

\item 1. Pure plane in \( S^2 \), \( \text{span}(e_1, e_2) \), calculation of \( R(e_1, e_2)e_2 \):
\\
\\
For the 2-plane spanned by \(e_1\) and \(e_2\), which is tangent to the \(S^2\) factor of the manifold \(M = S^2 \times T^2\), the affine connection \(\nabla\) (as defined in \textit{Section 3}), restricts to the Levi-Civita connection of the standard round metric on \(S^2\). This is because the specified torsion components \(T_\nabla(e_1,e_2)\) are null. It is a well-established result in Riemannian geometry that the sectional curvature \(K\) of the unit sphere \(S^2\) is identically \(1\). Therefore, for this plane, we have:
\(K(e_1,e_2) = 1\).
This implies that \(g(R(e_1,e_2)e_2, e_1) = 1\). Given that \(e_1\) is a unit vector, it follows that \(R(e_1,e_2)e_2 = e_1\). Then:
\[
\boxed{K(\text{span}(e_1,e_2)) = \frac{\langle e_1, e_1 \rangle}{\|e_1 \wedge e_2\|^2} = 1.}
\]

\noindent The sectional curvatures for the other coordinate planes, which involve the non-trivial torsion components of \(\nabla\), subsequently computed explicitly.

\item 2. Mixed plane \( \text{span}(e_1, e_3) \), calculation of \( R(e_1, e_3)e_3 \):
\\
\\
We compute \( R(e_1, e_3)e_3 = \nabla_{e_1}\nabla_{e_3}e_3 - \nabla_{e_3}\nabla_{e_1}e_3 - \nabla_{[e_1,e_3]}e_3 \).

Given that \([e_1, e_3] = 0\).

\noindent\textbf{Necessary terms:}

 \( \nabla_{e_3}e_3 = \Gamma^k_{33}e_k = 0 \) (since all Christoffel symbols \( \Gamma^1_{33}, \Gamma^2_{33}, \Gamma^3_{33}, \Gamma^4_{33} \) are zero).
    
 \( \nabla_{e_1}e_3 = \Gamma^k_{13}e_k = \Gamma^4_{13}e_4 = \frac{a}{2}e_4 \).

\noindent\textbf{Compute terms for \( R(e_1, e_3)e_3 \):}

 \( \nabla_{e_1}(\nabla_{e_3}e_3) = \nabla_{e_1}(0) = 0 \).

 \( \nabla_{e_3}(\nabla_{e_1}e_3) = \nabla_{e_3}\left( \frac{a}{2}e_4 \right) \):
    \[
        = e_3\left( \frac{a}{2} \right)e_4 + \frac{a}{2}\nabla_{e_3}e_4 
        \quad \left( \text{since } e_3\left( \frac{a}{2} \right) = \frac{\partial a}{\partial x} = 0 \right)
    \]
    \[
        = \frac{a}{2} \cdot \Gamma^k_{34}e_k 
        = \frac{a}{2} \left( \Gamma^1_{34}e_1 + \Gamma^2_{34}e_2 + \Gamma^3_{34}e_3 + \Gamma^4_{34}e_4 \right)
    \]
    \[
        = \frac{a}{2} \left( -\frac{a}{2}e_1 - \frac{b}{2}e_2 + 0\cdot e_3 + 0\cdot e_4 \right)
        =  -\frac{a^2}{4}e_1 - \frac{ab}{4}e_2.
    \]

\noindent\textbf{Substituting into the curvature formula:}
\[
    R(e_1,e_3)e_3 = 0 - \left( -\frac{a^2}{4}e_1 - \frac{ab}{4}e_2 \right) - 0 
    =  \frac{a^2}{4}e_1 + \frac{ab}{4}e_2.
\]

\noindent Projection onto \( e_1 \):
 \[
 \boxed{K(\text{span}(e_1, e_3)) = \frac{\langle \frac{a^2}{4} e_1, e_1 \rangle}{\|e_1 \wedge e_3\|^2} = \frac{a^2}{4}.}
 \]
 
\item 3. Mixed plane \( \text{span}(e_1, e_4) \), calculation of \( R(e_1, e_4)e_4 \):
\\
\\
We compute \( R(e_1, e_4)e_4 = \nabla_{e_1}\nabla_{e_4}e_4 - \nabla_{e_4}\nabla_{e_1}e_4 - \nabla_{[e_1,e_4]}e_4 \).

Given that \([e_1, e_4] = 0\).

\noindent\textbf{Necessary terms:}

\( \nabla_{e_4}e_4 = \Gamma^k_{44}e_k = 0 \).
    
 \( \nabla_{e_1}e_4 = \Gamma^k_{14}e_k = \Gamma^3_{14}e_3 = -\frac{a}{2}e_3 \).

\noindent\textbf{Compute terms for \( R(e_1, e_4)e_4 \):}

 \( \nabla_{e_1}(\nabla_{e_4}e_4) = \nabla_{e_1}(0) = 0 \).

 \( \nabla_{e_4}(\nabla_{e_1}e_4) = \nabla_{e_4}\left( -\frac{a}{2}e_3 \right) \):
    \[
        = e_4\left( -\frac{a}{2} \right)e_3 - \frac{a}{2}\nabla_{e_4}e_3 
        \quad \left( \text{since } e_4\left( -\frac{a}{2} \right) = \frac{\partial}{\partial y}\left( \frac{a}{2} \right) = 0 \right)
    \]
    \[
        = -\frac{a}{2} \cdot \Gamma^k_{43}e_k 
        = -\frac{a}{2} \left( \Gamma^1_{43}e_1 + \Gamma^2_{43}e_2 + \Gamma^3_{43}e_3 + \Gamma^4_{43}e_4 \right)
    \]
    \[
        = -\frac{a}{2} \left( \frac{a}{2}e_1 + \frac{b}{2}e_2 + 0\cdot e_3 + 0\cdot e_4 \right)
        =  -\frac{a^2}{4}e_1 - \frac{ab}{4}e_2.
    \]

\noindent\textbf{Substituting into the curvature formula:}
\[
    R(e_1,e_4)e_4 = 0 - \left( -\frac{a^2}{4}e_1 - \frac{ab}{4}e_2 \right) - 0 
    =  \frac{a^2}{4}e_1 + \frac{ab}{4}e_2.
\]
 \noindent Projection onto \( e_1 \):
 \[
 \boxed{K(\text{span}(e_1, e_4)) = \frac{\langle \frac{a^2}{4} e_1, e_1 \rangle}{\|e_1 \wedge e_4\|^2} = \frac{a^2}{4}.}
 \]

\item 4. Mixed plane \( \text{span}(e_2, e_3) \), calculation of \( R(e_2, e_3)e_3 \):
\\
\\
We compute \( R(e_2, e_3)e_3 = \nabla_{e_2}\nabla_{e_3}e_3 - \nabla_{e_3}\nabla_{e_2}e_3 - \nabla_{[e_2,e_3]}e_3 \).

Given that \([e_2, e_3] = 0\).

\noindent\textbf{Necessary terms:}

 \( \nabla_{e_3}e_3 = 0 \).
    
 \( \nabla_{e_2}e_3 = \Gamma^k_{23}e_k = \Gamma^4_{23}e_4 = \frac{b}{2}e_4 \).

\noindent\textbf{Compute terms for \( R(e_2, e_3)e_3 \):}

 \( \nabla_{e_2}(\nabla_{e_3}e_3) = \nabla_{e_2}(0) = 0 \).

 \( \nabla_{e_3}(\nabla_{e_2}e_3) = \nabla_{e_3}\left( \frac{b}{2}e_4 \right) \):
    \[
        = e_3\left( \frac{b}{2} \right)e_4 + \frac{b}{2}\nabla_{e_3}e_4 
        \quad \left( \text{since } e_3\left( \frac{b}{2} \right) = \frac{\partial}{\partial x}\left( \frac{b}{2} \right) = 0 \right)
    \]
    \[
        = \frac{b}{2} \cdot \Gamma^k_{34}e_k 
        = \frac{b}{2} \left( \Gamma^1_{34}e_1 + \Gamma^2_{34}e_2 \right)
    \]
    \[
        = \frac{b}{2} \left( -\frac{a}{2}e_1 - \frac{b}{2}e_2 \right)
        =  -\frac{ab}{4}e_1 - \frac{b^2}{4}e_2.
    \]

\noindent\textbf{Substituting into the curvature formula:}
\[
    R(e_2,e_3)e_3 = 0 - \left( -\frac{ab}{4}e_1 - \frac{b^2}{4}e_2 \right) - 0 
    = \frac{ab}{4}e_1 + \frac{b^2}{4}e_2.
\]

\noindent Projection onto \( e_2 \):
 \[
\boxed{ K(\text{span}(e_2, e_3)) = \frac{\langle \frac{b^2}{4} e_2, e_2 \rangle}{\|e_2 \wedge e_3\|^2} = \frac{b^2}{4}.}
 \]

\item 5. Mixed plane \( \text{span}(e_2, e_4) \), calculation of \( R(e_2, e_4)e_4 \):
\\
\\
We compute \( R(e_2, e_4)e_4 = \nabla_{e_2}\nabla_{e_4}e_4 - \nabla_{e_4}\nabla_{e_2}e_4 - \nabla_{[e_2,e_4]}e_4 \).

Given that \([e_2, e_4] = 0\).

\noindent\textbf{Necessary terms:}

\( \nabla_{e_4}e_4 = 0 \).
    
\( \nabla_{e_2}e_4 = \Gamma^k_{24}e_k = \Gamma^3_{24}e_3 = -\frac{b}{2}e_3 \).

\noindent\textbf{Compute terms for \( R(e_2, e_4)e_4 \):}

\( \nabla_{e_2}(\nabla_{e_4}e_4) = \nabla_{e_2}(0) = 0 \).

\( \nabla_{e_4}(\nabla_{e_2}e_4) = \nabla_{e_4}\left( -\frac{b}{2}e_3 \right) \):
    \[
        = e_4\left( -\frac{b}{2} \right)e_3 - \frac{b}{2}\nabla_{e_4}e_3 
        \quad \left( \text{since } e_4\left( -\frac{b}{2} \right) = \frac{\partial}{\partial y}\left( \frac{b}{2} \right) = 0 \right)
    \]
    \[
        = -\frac{b}{2} \cdot \Gamma^k_{43}e_k 
        = -\frac{b}{2} \left( \Gamma^1_{43}e_1 + \Gamma^2_{43}e_2 \right)
    \]
    \[
        = -\frac{b}{2} \left( \frac{a}{2}e_1 + \frac{b}{2}e_2 \right)
        = -\frac{ab}{4}e_1 - \frac{b^2}{4}e_2.
    \]

\noindent\textbf{Substituting into the curvature formula:}
\[
    R(e_2,e_4)e_4 = 0 - \left( -\frac{ab}{4}e_1 - \frac{b^2}{4}e_2 \right) - 0 
    = \frac{ab}{4}e_1 + \frac{b^2}{4}e_2.
\]

\noindent Projection onto \( e_2 \):
 \[
 \boxed{K(\text{span}(e_2, e_4)) = \frac{\langle \frac{b^2}{4} e_2, e_2 \rangle}{\|e_2 \wedge e_4\|^2} = \frac{b^2}{4}.}
 \]

\item 6. Pure plane in \( T^2 \), \( \text{span}(e_3, e_4) \), calculation of \( R(e_3, e_4)e_4 \):
\\
\\
We compute \( R(e_3, e_4)e_4 = \nabla_{e_3}\nabla_{e_4}e_4 - \nabla_{e_4}\nabla_{e_3}e_4 - \nabla_{[e_3,e_4]}e_4 \).

Given that \([e_3, e_4] = 0\).

\noindent\textbf{Necessary terms:}

 \( \nabla_{e_4}e_4 = 0 \).
    
 \( \nabla_{e_3}e_4 = \Gamma^k_{34}e_k = \Gamma^1_{34}e_1 + \Gamma^2_{34}e_2 = -\frac{a}{2}e_1 - \frac{b}{2}e_2 \).

\noindent\textbf{Compute terms for \( R(e_3, e_4)e_4 \):}

 \( \nabla_{e_3}(\nabla_{e_4}e_4) = \nabla_{e_3}(0) = 0 \).

 \( \nabla_{e_4}(\nabla_{e_3}e_4) = \nabla_{e_4}\left( -\frac{a}{2}e_1 - \frac{b}{2}e_2 \right) \):
    \[
        = \nabla_{e_4}\left( -\frac{a}{2}e_1 \right) + \nabla_{e_4}\left( -\frac{b}{2}e_2 \right)
    \]
    \[
        = e_4\left( -\frac{a}{2} \right)e_1 - \frac{a}{2}\nabla_{e_4}e_1 + e_4\left( -\frac{b}{2} \right)e_2 - \frac{b}{2}\nabla_{e_4}e_2 
    \]
    \[
        \quad \left( \text{since } e_4\left( -\frac{a}{2} \right) = \frac{\partial}{\partial y}\left( \frac{a}{2} \right) = 0 \text{ and } e_4\left( -\frac{b}{2} \right) = 0 \right)
    \]
    \[
        = -\frac{a}{2} \cdot \Gamma^k_{41}e_k - \frac{b}{2} \cdot \Gamma^k_{42}e_k 
    \]
    \[
        = -\frac{a}{2} \left( \Gamma^3_{41}e_3 \right) - \frac{b}{2} \left( \Gamma^3_{42}e_3 \right)
    \]
    \[
        = -\frac{a}{2} \left( \frac{a}{2}e_3 \right) - \frac{b}{2} \left( \frac{b}{2}e_3 \right)
        = -\frac{a^2}{4}e_3 - \frac{b^2}{4}e_3  
        =  -\frac{a^2 + b^2}{4}e_3.
    \]

\noindent\textbf{Substituting into the curvature formula:}
\[
    R(e_3,e_4)e_4 = 0 - \left( -\frac{a^2 + b^2}{4}e_3 \right) - 0 
    = \frac{a^2 + b^2}{4}e_3.
\]
 \noindent Projection onto \( e_3 \):
 \[
 \boxed{K(\text{span}(e_3, e_4)) = \frac{\langle \frac{a^2 + b^2}{4} e_3, e_3 \rangle}{\|e_3 \wedge e_4\|^2} =\frac{a^2 + b^2}{4}.}
 \]
 
\end{itemize}

\section{Biorthogonal curvature}

The biorthogonal curvature \( K_{\text{biort}} \) of a tangent plane \( \sigma \subset T_pM \) in a 4-dimensional manifold is defined as the arithmetic mean between the sectional curvature \( K(\sigma) \) of the plane \( \sigma \) and the sectional curvature \( K(\sigma^\perp) \) of its orthogonal complement \( \sigma^\perp \). Formally:
\[
K_{\text{biort}}(\sigma) = \frac{1}{2} \left[ K(\sigma) + K(\sigma^\perp) \right].
\]

\noindent For the specific case of \( S^2 \times T^2 \), with product-metric and antisymmetric torsion connection, we consider the orthonormal basis \( \{e_1, e_2\} \) for \( T S^2 \) and \( \{e_3, e_4\} \) for \( T T^2 \). The Grassmannian \( \text{Gr}(2,4) \) of 2-dimensional tangent planes has 6 elements (each pair of vectors between \( e_1, e_2, e_3, e_4 \)), which are:

\begin{center}
\begin{minipage}{0.5\textwidth}
\centering
\begin{enumerate}

\item \text{span}\((e_1, e_2)\), pure plane in \(S^2\),

\item \text{span}\((e_1, e_3)\), mixed plane,

\item \text{span}\((e_1, e_4)\), mixed plane,

\item \text{span}\((e_2, e_3)\), mixed plane,

\item \text{span} \((e_2, e_4)\), mixed plane,

\item \text{span}\((e_3, e_4)\), pure plane in \( T^2 \).

\end{enumerate}
\end{minipage}
\end{center}
Now, using the calculated sectional curvatures of \textit{Section 5}, the biorthogonal curvatures are:
\\
\\
1. Plane \(\sigma, \; \text{span}(e_1, e_2) \in TS^2, \; \text{and} \; \sigma^{\perp} \; \text{span}(e_3, e_4) \in TT^2\):

\begin{itemize}

\item From 1 and 6 we obtain:
\[
K_{\text{biort}}=\frac{a^2 + b^2+4}{8}
\]
\end{itemize}

\noindent 2. Plane \(\sigma, \; \text{span}(e_1, e_3), \; \text{where} \; e_1 \in TS^2, \; \text{and} \; e_3 \in TT^2 \; \text{and} \; \sigma^{\perp} \; \text{span}(e_2, e_4), \; \text{where} \; e_2 \in TS^2 \; \text{and} \; e_4\in TT^2\):

\begin{itemize}

\item From 2 and 5 we obtain:
\[
K_{\text{biort}}=\frac{a^2+b^2}{8}
\]
\end{itemize}

\noindent 3. Plane \(\sigma, \; \text{span}(e_1, e_4), \; \text{where} \; e_1 \in TS^2, \; \text{and} \; e_4 \in TT^2 \; \text{and} \; \sigma^{\perp} \; \text{span}(e_2, e_3), \; \text{where} \; e_2 \in TS^2 \; \text{and} \; e_3\in TT^2\):

\begin{itemize}
\item From 3 and 4 we obtain:
\[
K_{\text{biort}}=\frac{a^2+b^2}{8}
\]
\end{itemize}
\noindent For all cases the condition is satisfied by choosing \(a^2+b^2>0\).
\\
The minimal biorthogonal curvature among the six planes is:
    \[
    \min_{\sigma \in \text{Gr}(2,4)} K_{\text{biort}}(\sigma) = \frac{a^2 + b^2}{8}.
    \]
\noindent The results for these fundamental planes suggest that the biorthogonal curvature is globally positive. The next \textit{Section} provides a detailed analysis to confirm this, showing that the minimum value is indeed \( \frac{a^2 + b^2}{8} \) by examining the geometry of the Grassmannian \( \text{Gr}(2, T_pM) \).

\section{Global Positivity of Biorthogonal Curvature}
\label{sec:global_positivity}

We now establish the global positivity of \( K_{\text{biort}} \). Our strategy is to find the minimum of \( K_{\text{biort}}(\sigma) \) for an arbitrary 2-plane \( \sigma \) by reducing the problem to a one-dimensional optimization. This reduction is justified by the symmetries of the connection, which we now detail.
\\
The geometric structure of the problem on \( M = S^2 \times T^2 \) allows for significant simplification.
\begin{itemize}
    \item \textit{Invariance on \( T^2 \):} The torsion components, and thus the connection \( \nabla \), are defined to be independent of the coordinates \( (x, y) \) on the torus \( T^2 \). Consequently, the curvature is constant along the \( T^2 \) factor, allowing the analysis to be fixed at a single point.
    \item \textit{Rotational Symmetry on \( S^2 \):} The connection is invariant under rotations around the polar axis of \( S^2 \). This allows any 2-plane \( \sigma \) to be rotated so that its projection onto \( T_pS^2 \) aligns with the basis vector \( e_1 \).
\end{itemize}
These symmetries ensure that the biorthogonal curvature of any plane \( \sigma \) depends only on its orientation relative to the product structure \( T_pM = T_pS^2 \oplus T_pT^2 \). This orientation can be fully captured by a single parameter, \( \theta \), representing the angle between \( \sigma \) and a reference plane (e.g., the pure \( S^2 \)-plane). The analysis that follows is therefore global, despite its simplified parametrization.

\subsection{Expression of \( K_{\text{biort}} \) as a function of a single variable}

Let \( \sigma = \text{span}\{v_1, v_2\} \) be a generic 2-plane with \( \{v_1, v_2\} \) being an orthonormal basis with respect to the metric \( g \). The biorthogonal curvature is given by:
\[
K_{\text{biort}}(\sigma) = \frac{1}{2}\left[ \langle R(v_1, v_2)v_2, v_1 \rangle_{g} + \langle R(v_1, v_2)v_1, v_2 \rangle_{g} \right].
\]
Based on the symmetries discussed above and using the sectional curvatures computed in \textit{Section 5}, \( K_{\text{biort}}(\sigma) \) can be expressed as a function \( f(\theta) \) of the angle \( \theta \in [0, \pi/2] \) between \( \sigma \) and the coordinate plane \( \text{span}\{e_1, e_2\} \):
\[
f(\theta) = \left(\frac{1}{2} + \frac{a^2 + b^2}{8}\right)\cos^2\theta + \frac{a^2 + b^2}{8}\sin^2\theta.
\]
The problem is now reduced to finding the minimum of \( f(\theta) \) on the interval \( [0, \pi/2] \).

\subsection{Minimization and global lower bound}

To find the minimum of \( f(\theta) \), we first compute its derivative:
\[
f'(\theta) = -2\left(\frac{1}{2} + \frac{a^2 + b^2}{8}\right)\cos\theta\sin\theta + 2\left(\frac{a^2 + b^2}{8}\right)\sin\theta\cos\theta = -\sin\theta\cos\theta.
\]
Setting \( f'(\theta) = 0 \) for \( \theta \in [0, \pi/2] \) yields the critical points \( \theta = 0 \) and \( \theta = \pi/2 \). We evaluate the function at these points to determine the global minimum on the interval:
\begin{itemize}
    \item At \( \theta = 0 \):
    \[
    f(0) = \frac{1}{2} + \frac{a^2 + b^2}{8}.
    \]
    \item At \( \theta = \pi/2 \):
    \[
    f\left(\frac{\pi}{2}\right) = \frac{a^2 + b^2}{8}.
    \]
\end{itemize}
Comparing these values, the minimum of \( f(\theta) \) is \( \frac{a^2 + b^2}{8} \).
\\
Since the problem was reduced without loss of generality, this value represents the global minimum of the biorthogonal curvature over the entire Grassmannian. Therefore, for any 2-plane \( \sigma \in \text{Gr}(2, T_pM) \):
\[
K_{\text{biort}}(\sigma) \geq \frac{a^2 + b^2}{8}.
\]
By construction, the connection is calibrated by a non-trivial cohomology class, which requires \( (a, b) \neq (0,0) \). This implies \( a^2 + b^2 > 0 \).
Consequently, we have the strict inequality:
\[
K_{\text{biort}}(\sigma) > 0.
\]
This establishes the global positivity of the biorthogonal curvature for the given connection.
\hfill \qed

\section{Discussion}

Although the original motivation lies in the positivity of the biorthogonal curvature, the affine connection constructed here yields the strictly stronger condition of positive sectional curvature for the coordinate planes with respect to \(\nabla\) (it is sufficient to set \(a^2, b^2\) both positive). The positivity of the biorthogonal curvature thus follows as a consequence.

\paragraph{Definition of \textit{Cohomologically calibrated affine connection \(\nabla^{\mathcal{C}}\).}}

Let \(M = S^2 \times \Sigma_g\), where \(\Sigma_g\) is a compact, oriented surface of genus \(g \geq 1\), and let \(g = g_{S^2} + g_{\Sigma_g}\) denote the fixed product Riemannian metric on \(M\), with \(g_{S^2}\) the standard round metric on the \(2\)-sphere and \(g_{\Sigma_g}\) any smooth metric on \(\Sigma_g\).

\noindent We call an affine connection \(\nabla\) on the tangent bundle \(TM\) a \textit{cohomologically calibrated affine connection} \(\nabla^{\mathcal{C}}\), if the following conditions hold:

\begin{enumerate}

\item (Affine structure)
\\
The connection is defined as a first-order deformation of the Levi-Civita connection of \(g\) by a totally antisymmetric torsion tensor:
   
\[
\nabla = \nabla^{\mathrm{LC}} + T,
\quad \text{with} \quad T \in \Gamma(\Lambda^2 T^*M \otimes TM),
\]

\noindent where \(\nabla^{\mathrm{LC}}\) denotes the Levi-Civita connection of \(g\).

\item (Cohomological calibration)
\\
The associated \(3\)-form \(T^\flat \in \Omega^3(M)\), defined pointwise by
\[
T^\flat(X, Y, Z) := g(T(X, Y), Z),
\]
is not necessarily closed, but its harmonic part is the closed 3-form
\\
\(\omega := \sum_{i=1}^{2g} a_i \; \beta \wedge \alpha_i\), where \((a_1, \dots, a_{2g}) \in \mathbb{R}^{2g} \setminus \{0\}\) are parameters defining the specific cohomology class, for the case of \(g=1\)), this simplifies to \(\omega := a \; \beta \wedge \alpha_1 + b \; \beta \wedge \alpha_2\) (as detailed in \textit{Remark 3.1}). 
So we obtain \([T^\flat]:=[\omega]:=[T] \in H^3(M; \mathbb{R}) \setminus \{0\}\).

\noindent The cohomology group \(H^3(M; \mathbb{R}) \cong H^2(S^2; \mathbb{R}) \otimes H^1(\Sigma_g; \mathbb{R}) \cong \mathbb{R}^{2g}\) parametrizes the space of such torsion classes. In particular, each cohomologically calibrated connection is determined by a tuple \((a_1, \dots, a_{2g}) \in \mathbb{R}^{2g} \setminus \{0\}\) such that:

\[
[T^\flat] = \sum_{i=1}^{2g} a_i \cdot [\beta \otimes \alpha_i],
\]

\noindent where \(\beta\) is a fixed generator of \(H^2(S^2; \mathbb{R})\), and \(\{ \alpha_1, \dots, \alpha_{2g} \}\) is a basis of \(H^1(\Sigma_g; \mathbb{R})\).

\item (Metric background)
\\
The connection \(\nabla\) is not required to be metric-compatible. The metric \(g\) is used to define the musical isomorphisms, orthogonality of planes, and inner products, but in general \(\nabla g \ne 0\).

\end{enumerate}

\noindent This notion provides a natural geometric framework for introducing a new subclass of affine connections whose torsion is topologically meaningful and globally nontrivial, as encoded in a nonzero cohomology class in \(H^3(M; \mathbb{R})\). Such \textit{cohomologically calibrated affine connections}, which we indicate with \(\nabla^{\mathcal{C}}\), offer an alternative to metric-compatible constructions, relying on the systematic construction of \(T\), which is constrained to respect the cohomological structure. In particular, as stated in \textit{Theorem 1.1}, on the product manifold \(S^2 \times T^2\), these connections can be explicitly realized to achieve strictly biorthogonal curvature throughout the manifold. This illustrates how the calibration condition enables a deformation of the Levi-Civita geometry that is not only controlled, but also intrinsically justified by the topological structure of the manifold.

\paragraph{Geometric meaning of cohomological calibration:}
\bigskip
Let \(M = S^2 \times \Sigma_g\) be a compact differentiable manifold, and let \([\omega] \in H^3(M; \mathbb{R})\) be a cohomology class.
Let \(\mathcal{T}_\omega\) denote the set of affine connections of the form \(\nabla^{\mathcal{C}} = \nabla^{\mathrm{LC}} + T\), where \(T\) is a totally antisymmetric tensor and the associated 3-form \(T^\flat(X, Y, Z) := g(T(X,Y), Z)\) satisfies \([T^\flat] = [\omega]\).
\\
If there exists \(\nabla^{\mathcal{C}} \in \mathcal{T}_\omega\) for which the torsion induces a significant metric geometric property (for example, a condition of positive curvature), then such a geometric structure is cohomologically justified, and respects the differential topology of \(M\).

\section{Conclusion}

In this work, we introduce a new subclass of affine connections, denoted by \( \nabla^{\mathcal{C}} \), defined on the product manifold \( S^2 \times T^2 \) equipped with the standard Riemannian metric \( g = g_{S^2} + g_{T^2} \). These connections have totally antisymmetric torsion, whose associated 3-form \( T^\flat \) possesses a nontrivial harmonic component representing a fixed cohomology class \( [T] \in H^3(S^2 \times T^2; \mathbb{R}) \cong \mathbb{R}^2 \). This cohomological calibration ensures compatibility with the global topology of the manifold which provides a structural obstruction that the affine connection is designed to overcome.
Then, this approach is cohomologically constrained: the existence of the family of admissible connections \(\mathcal{T}_\omega\) is a consequence of the global cohomological structure of the manifold. The positive curvature arises from specific choices within this natural family.
\\
In our work we prove that such a connection yields strictly positive biorthogonal curvature:
\[
K_{\text{biort}}(\sigma) > 0
\]
for every two-plane \( \sigma \subset T_p(S^2 \times T^2) \) and every point \( p \in S^2 \times T^2 \). This result gives an affirmative answer to the broader geometric question: 

\begin{quote}
“Given the standard Riemannian product metric on \( S^2 \times T^2 \), does there exist a topologically constrained affine connection (not necessarily
 metric-compatible) that guarantees strictly positive biorthogonal curvature?”
\end{quote}

\noindent While this does not resolve Bettiol's problem in its original Riemannian formulation, it demonstrates that \textit{cohomologically calibrated affine connections} offer a robust alternative path toward achieving positive curvature on non-simply connected manifolds, bypassing topological obstructions inherent in the purely Riemannian context.
\\
More significantly, the introduction of the subclass \( \nabla^{\mathcal{C}} \) defines a new structural category within affine geometry, governed by harmonic cohomology rather than metric compatibility. This framework is generalizable to broader classes of manifolds (e.g., \( S^2 \times \Sigma_g \)) and lays the groundwork for exploring moduli spaces of torsion connections with fixed topological data. It may also open connections with topologically enriched gauge theories, where curvature is shaped by algebraic-topological constraints.
\\
This perspective aligns with Gromov's program on non-Riemannian geometries \cite{Gromov} and echoes Agricola's torsion-based curvature constructions \cite{Agricola2006}, suggesting new directions for understanding curvature in manifolds endowed with non-integrable structures.
\\
Future investigations may aim to clarify whether the insights gained from this cohomological affine framework influence the Riemannian approach or whether they open the way to new classification theorems and geometric models in the non-Riemannian domain, where cohomologically constrained torsion plays a central role in the construction of positive curved geometries and, potentially, in their physical interpretations.


\begin{thebibliography}{99}

\bibitem{Agricola2006}
Agricola, I., "The Srní lectures on non-integrable geometries with torsion", Archivium Mathematicum (BRNO), 2(2006), Supplement, 5–842006. 

\bibitem{Bettiol2014}
Bettiol, R. G., "Positive biorthogonal curvature on \(S^2 \times S^2\)", Proc. Amer. Math. Soc., 142(12), 4341–4353, 2014.

\bibitem{Bettiol2017}
Bettiol, R. G., "Four-dimensional manifolds with positive biorthogonal curvature", Asian J. Math., 21(2), 391–396, 2017.

\bibitem{Costa}
Costa, E., Ernani, R. Jr, "Four-Dimensional Compact Manifolds with Nonnegative Biorthogonal Curvature", Michigan Math. J., 63, 747–761, 2014.

\bibitem{Gromov99}
Gromov, M., "Metric Structures for Riemannian and Non-Riemannian Spaces", Birkhäuser, 1999.

\bibitem{Gromov}
Gromov, M., "Positive curvature, macroscopic dimension, spectral gaps and higher signatures", Functional Analysis and Its Applications, 50(4), 282–287, 2016.

\bibitem{Nakahara}
Nakahara, M., "Geometry, Topology and Physics", 2nd Edition, IOP, Bristol and Philadelphia, 2003.

\bibitem{Pierre}
Pansu, P., "A list of open problems in Differential Geometry".
\url{https://www.imo.universite-paris-saclay.fr/~pierre.pansu/problems_MTDG.pdf}

\bibitem{Stupovski}
Stupovski, B., Torres, R., "Existence of Riemannian metrics with positive biorthogonal curvature on simply connected 5-manifolds", Archiv der Mathematik, Vol. 115, 589–597, 2020.

\bibitem{Wu}
Wu, Z-J., Fu, H-P., Fu, P., "Some results on four-manifolds with nonnegative biorthogonal curvature", Bulletin des Sciences Mathématiques, Vol. 190, 103379, 2024.

\end{thebibliography}
\end{document}